\documentclass[journal]{IEEEtran}
\usepackage[utf8]{inputenc}

\title{Optimal Carbon Taxes for Emissions Targets in the Electricity Sector}
\author{Daniel~J.~Olsen,~\IEEEmembership{Student~Member,~IEEE,}
        Yury~Dvorkin,~\IEEEmembership{Member,~IEEE,}
        Ricardo~Fern\'andez-Blanco,
        and Miguel~A.~Ortega-Vazquez,~\IEEEmembership{Senior Member,~IEEE}}
\date{September 2017}

\usepackage{amsmath,amssymb,amsthm}     
\usepackage{cite}
\usepackage{xcolor}                     
\usepackage{graphicx}
\usepackage{float}
\usepackage{array}                      
\usepackage{wasysym}                    

\usepackage{tikz}
\usepackage{textcomp}
\usepackage{hyperref}
\usepackage{lipsum}

\newcommand{\copyrighttext}{%
	\footnotesize \textcopyright 2018 IEEE. Personal use of this material is permitted. Permission from IEEE must be obtained for all other uses, in any current or future media, including reprinting/republishing this material for advertising or promotional purposes, creating new collective works, for resale or redistribution to servers or lists, or reuse of any copyrighted component of this work in other works.	DOI: \href{https://ieeexplore.ieee.org/document/8338157/}{10.1109/TPWRS.2018.2827333}, IEEE Transactions on Power Systems.}

\newcommand{\copyrightnotice}{%
	\begin{tikzpicture}[remember picture,overlay]
	\node[anchor=south,yshift=8pt] at (current page.south) {\fbox{\parbox{\dimexpr\textwidth-\fboxsep-\fboxrule\relax}{\copyrighttext}}};
	\end{tikzpicture}%
}

\newenvironment{ldescription}[1]        
  {\begin{list}{}%
   {\renewcommand\makelabel[1]{##1\hfill}%
   \settowidth\labelwidth{\makelabel{#1}}%
   \setlength\leftmargin{\labelwidth}
   \addtolength\leftmargin{\labelsep}}}
  {\end{list}}

\begin{document}

\maketitle
\copyrightnotice

\begin{abstract}
    The most dangerous effects of anthropogenic climate change can be mitigated by using emissions taxes or other regulatory interventions to reduce greenhouse gas (GHG) emissions.  This paper takes a regulatory viewpoint and describes the Weighted Sum Bisection method to determine the lowest emission tax rate that can reduce the anticipated emissions of the power sector below a prescribed, regulatorily-defined target.  This bi-level method accounts for a variety of operating conditions via stochastic programming and remains computationally tractable for realistically large planning test systems, even when binary commitment decisions and multi-period constraints on conventional generators are considered.
    
    Case studies on a modified ISO New England test system demonstrate that this method reliably finds the minimum tax rate that meets emissions targets.  In addition, it investigates the relationship between system investments and the tax-setting process.  Introducing GHG emissions taxes increases the value proposition for investment in new cleaner generation, transmission, and energy efficiency; conversely, investing in these technologies reduces the tax rate required to reach a given emissions target.
\end{abstract}


\vspace{-4pt}

\section*{Nomenclature}

\vspace{-8pt}
\subsection*{Sets and Indices}
\begin{ldescription}{$xxxxx$}
    \item [$A$] Set of representative days, indexed by $a$.
    \item [$B$] Set of transmission network buses, indexed by $b$.
    \item [$I$] Set of generating units, indexed by $i$.
    \item [$L$] Set of transmission lines, indexed by $l$.
    \item [$R$] Subset of renewable generators ($R \subset I$).
    \item [$S$] Set of generator power output blocks, indexed by $s$.
    \item [$T$] Set of time intervals, indexed by $t$ or $\tau$.
\end{ldescription}

\vspace{-14pt}
\subsection*{Parameters}
\begin{ldescription}{$xxxxx$}
    \item [$b_{i,s}$] Marginal cost of block $s$ of generator $i$ (\$/MWh).
    \item [$C_{i}^{\text{min}}$] Minimum cost of generator $i$ (\$/h).
    \item [$C_{i}^{\text{su}}$] Start-up cost of generator $i$ (\$).
    \item [$d_{b,t,a}$] Demand at bus $b$, time $t$, day $a$ (MW).
    \item [$d_{t,a}^{\text{ramp}}$] Load ramp requirement at time $t$, day $a$ (MW/h).
    \item [${E}^{\text{max}}$] Regulator's GHG emission target (tons).
    \item [$E_{i}^{\text{min}}$] Minimum GHG emissions of generator $i$ (tons/h).
    \item [$E_{i}^{\text{su}}$] Start-up GHG emissions of generator $i$ (tons).
    \item [$f_{l}^{\text{max}}$] Capacity of transmission line $l$ (MW).
    \item [$g_{i}^{\text{max}}$] Maximum power output of generator $i$ (MW).
    \item [$g_{i}^{\text{min}}$] Minimum power output of generator $i$ (MW).
    \item [$g_{i,s}^{\text{max}}$] Maximum power output of block $s$, generator $i$ (MW).
    \item [$g_{i}^{\text{down}}$] Minimum down-time of generator $i$ (h).
    \item [$g_{i}^{\text{up}}$] Minimum up-time of generator $i$ (h).
    \item [$h_{i,s}$] Marginal GHG emissions of block $s$, generator $i$ (tons/MWh).
    \item [$m_{l,b}^{\text{line}}$] Line connection map. ${m}_{lb}^{\text{line}} = 1$ if line $l$ starts at bus $b$, $= -1$ if line $l$ ends in bus $b$, $0$ otherwise.
    \item [$m_{i,b}^{\text{unit}}$] Unit map. ${m}_{i,b}^{\text{unit}} = 1$ if generator $i$ is located at bus $b$, $0$ otherwise.
    \item [$P^{\text{CO}_2}$] GHG emissions tax rate (\$/ton-$\text{CO}_{2}$e).
    \item [$P^{\text{load}}$] Load shed penalty (\$/MWh).
    \item [$P^{\text{ren}}$] Renewable generation shed penalty (\$/MWh).
    \item [$r_{i}^{\text{down}}$] Maximum down-ramp rate of generator $i$ (MW/h).
    \item [$r_{i}^{\text{up}}$] Maximum up-ramp rate of generator $i$ (MW/h).
    \item [$w_{t,a}^{\text{down}}$] Wind down-ramp requirements at time $t$, day $a$ (MW/h).
    \item [$w_{t,a}^{\text{up}}$] Wind up-ramp requirements at time $t$, day $a$ (MW/h).
    \item [$x_{l}$] Reactance of line $l$ ($\Omega$).
    \item [$\pi_{a}$] Probability of day $a$.
\end{ldescription}

\vspace{-14pt}
\subsection*{Variables}
\begin{ldescription}{$xxxxx$}
    \item [$C^{\text{gen}}$] System operator's generation cost (\$).
    \item [$C^{\text{shed}}$] System operator's shed cost (\$).
    \item [$E_{a}$] GHG Emissions for day $a$ (tons).
    \item [$E^{\text{total}}$] Total GHG emissions (tons).
    \item [$f_{l,t,a}$] Power flow on line $l$, time $t$, day $a$ (MW).
    \item [$g_{i,t,a}$] Power output of generator $i$, time $t$, day $a$ (MW).
    \item [$g_{i,s,t,a}$] Power output of generator $i$, block $s$, time $t$, day $a$ (MW).
    \item [$s_{b,t,a}^{\text{load}}$] Load shed at bus $b$, time $t$, day $a$ (MWh).
    \item [$s_{b,t,a}^{\text{ren}}$] Renewable generation shed at bus $b$, time $t$, day $a$ (MWh).
    \item [$u_{i,t,a}$] Binary variable for the commitment status of generator $i$, time $t$, day $a$.
    \item [$v_{i,t,a}$] Binary variable for the start-up of generator $i$, time $t$, day $a$.
    \item [$z_{i,t,a}$] Binary variable for the shut-down of generator $i$, time $t$, day $a$.
    \item [$\theta_{b,t,a}$] Voltage phase angle of bus $b$, time $t$, day $a$ (rad).
\end{ldescription}

\section{Introduction}

\subsection{Background}

The risks posed by anthropogenic climate change are dire, and organized effort is required in order to mitigate and eliminate, when possible, the effects \cite{ipcc_syr_2014}.  Pricing the emissions of greenhouse gases (GHGs) is a well-established approach to internalizing these negative externalities and should result in shifting the supply-demand equilibrium to a socially optimal point \cite{pigou}. Since the true costs from climate change are uncertain and hard to quantify with any precision (though attempts have been made, as in \cite{nas_social_cost}), one approach to create a price for emissions is to design policies that aim to reduce emissions to a level that is generally accepted to avoid the worst effects.  Unlike renewable portfolio standards or tax credits for renewable energy investment or production, this approach is directly targeted toward reducing GHG emissions.  As noted in \cite{deng2015}, subsidies for production of renewable energy can lead to negative bids by renewable generators, which may result in higher costs and emissions than if they bid zero-cost. According to a 2017 study by the World Bank \cite{worldbank2017}, there are 47 regional, national, and sub-national carbon pricing schemes implemented or scheduled for implementation, ranging from \$1-140/t$\text{CO}_{2}\text{e}$ and covering 15\% of global emissions.

The two main approaches to pricing emissions are a tax on GHG emissions (\textit{i.e.} a carbon tax) \cite{vehmas2005} and a cap-and-trade system \cite{egenhofer2007}.  A carbon tax sets a price directly with the goal of implicitly reducing emissions, while a cap-and-trade system sets emissions reductions explicitly, implicitly creating a price.  Each system has pros and cons: a cap-and-trade system can be more precise about the level of emissions reductions achieved, but requires complex rules regarding distributing, auctioning, and trading of allowances.  A carbon tax is simpler and may be easier to implement, but impact on emissions is less certain \cite{avi-yonah2009}, as the reactions to such a tax by the broader market (\textit{e.g.} generation and transmission investors, electricity consumers, and generation manufacturers) are difficult to model.
 
Secondary policy considerations are similar between the two: entities may purchase carbon offsets to reduce net emissions, tariffs can maintain competitiveness with jurisdictions without carbon pricing, and policies can be designed to be revenue-neutral.  Such policy design considerations for a carbon tax are discussed in \cite{metcalf2009,metcalf2017}. Additionally, pricing of carbon in either approach can lead to carbon `leakage', \textit{i.e.} the increased cost of producing goods in a jurisdiction with a carbon tax can lead to a shift in production toward jurisdictions with lower rates, or no tax at all \cite{sauma2012}.

The impact of carbon taxes on the economy and the environment have been widely studied, using various tax rates: the Brookings Institute in 2012 studied a tax which would begin at \$15/t$\text{CO}_{2}\text{e}$ with an annual escalator \cite{brookings2012}, the Congressional Budget Office in 2013 evaluated various rates between \$15-29/t$\text{CO}_{2}\text{e}$ \cite{cbo_tax}, and the Energy Information Agency in 2014 investigated rates of \$10 \& \$25/t$\text{CO}_{2}\text{e}$ \cite{aeo2014}.  However, these studies do not address at what rate carbon should be taxed.

This paper approaches the topic from a different angle.  Instead of studying the impact of a certain tax rate, we set a tax rate to achieve a certain environmental impact (\textit{e.g.} pledges from the Paris Agreement \cite{carbonbrief2017}) at minimal tax rate.  Specifically, we present an approach for setting the \emph{optimal carbon tax} for a given power system such that the resulting minimum-cost generator commitment and dispatch yields emissions that are at or below a specified target.  The consequence of a tax rate that is too low is failure to meet emissions targets, and the consequence of a rate too high is undue economic burden. Minimizing the tax rate also has practical motivations: lower tax rates are often more politically palatable (\textit{i.e.} more likely to be enacted, less likely to be repealed), and generally reduce rates of tax evasion \cite{liu2013}.

Some carbon taxing systems are designed to `recycle' the revenue received, either by investing in clean generation technologies or by reducing tax rates on other sectors of the economy to achieve revenue `neutrality'. However, care must be taken to account for the uncertainty in future carbon consumption, especially as carbon pricing tends to reduce consumption.

\subsection{Literature Survey}

Work including emissions into generator dispatch began in the 1970s with the concepts of minimum emissions dispatch \cite{gent1971}, pricing of emissions to include their impact in economic dispatch \cite{delson1974}, and varying emissions prices to investigate the tradeoffs between fuel costs and emissions \cite{zahavi1975}.  However, these studies focused on local effects of $\text{NO}_{x}$ and $\text{SO}_{x}$.

The concepts of `pseudo fuel prices' and algorithms for setting them are explored in \cite{lee1991}, and expanded in \cite{lee1992} to include periodic price adaptation in order to meet long-term fuel consumption targets as realizations differ from projections.  Similar algorithms are used to set weights based on emissions targets in \cite{ramanthan1994} for economic dispatch problems.  Several methods for coordinating long-term targets for emissions and short-term operations are discussed in \cite{gardner1995,fu2005}, but explicit emission pricing is absent, to the best of the authors' knowledge.

Explicit GHG pricing and its impact on optimal power flow problems are discussed in \cite{shao2010}. \cite{wei2014} presents a bi-level approach for setting a tax rate to achieve a GHG emissions target with minimal tax burden, but intertemporal constraints (\textit{e.g.} ramp rates, start-ups) are ignored. Without considering these constraints, the determined tax rate may not meet the desired target. `Optimal' tradeoffs between GHG emissions and cost according to a Nash bargaining process are developed in \cite{wei2015,wei2016}.

All of the above approaches contain deficiencies when it comes to setting a carbon tax rate with an eye on scheduling algorithms (\textit{i.e.} unit commitment models).  In short, one or more of the following is missing:
\begin{enumerate}
    \item Intertemporal variables and constraints (\textit{e.g.} ramp rate limitations, minimum up- and down-times).
    \item Explicit carbon pricing in dispatch/commitment.
    \item A method for setting an optimal carbon price.
\end{enumerate}

By contrast, in this work we propose a Weighted Sum Bisection (WSB) method, a computationally efficient approach, to set the minimal carbon tax rate that results in a power system meeting emissions targets while incorporating unit commitment, ramp rate limitations, and system flexibility and contingency reserve requirements.

\subsection{Contributions}

This work makes the following contributions:

\begin{enumerate}
    \item A bi-level planning model including unit commitment based on cyclic representative days, avoiding the need for assumptions about initial conditions.
    \item An efficient method for determining the minimal carbon tax rate which achieves emissions reductions targets.
    \item A demonstration of the computational efficiency of the proposed method and of the reciprocal relationship between tax rates for emissions targets and investment decisions.
\end{enumerate}

\section{Problem Formulation} \label{powersystem}

The problem is formulated as a bi-level planning problem, with the regulator's tax rate ($P^{\text{CO}_2}$) optimization in the upper level \eqref{upper_obj}-\eqref{te_max} and the system operator's unit commitment with carbon tax (UCCT), over a set of representative days, in the lower level \eqref{lower_obj}-\eqref{spill_limit}.  The UCCT includes the dc power flow approximation of the system power flows, penalties for shedding load and renewable generation, and reserve and ramping adequacy requirements. We assume an electricity market based on a unit commitment in which bids represent true fuel and tax costs.

\vspace{-20pt}
\begin{gather}
    \min P^{\text{CO}_2} \label{upper_obj} \\[-6pt]
    \text{subject to:} \hspace{100pt} \nonumber \\[-4pt]
    E^{\text{total}} \le E^{\text{max}} \label{te_max}
\end{gather}
\vspace{-21pt}
\begin{align}
    E^{\text{total}} \in \text{arg} \min \Big\{ C^{\text{shed}} + C^{\text{gen}} &+ P^{\text{CO}_2} E^{\text{total}} 
    \label{lower_obj} \\
    C^{\text{shed}} := \sum_{a \in A} \pi_a \sum_{t \in T} \sum_{b \in B}
    \Big( P^{\text{load}} s_{b,t,a}^{\text{load}} &+ \sum_{i \in R} P^{\text{ren}} s_{i,t,a}^{\text{ren}} \Big)
    \label{shed_def} \\
    C^{\text{gen}} := \sum_{a \in A} \pi_a \sum_{t \in T} \sum_{i \in I}
    \Big( C_{i}^{\text{min}} u_{i,t,a} &+ {C}_{i}^{\text{su}} v_{i,t,a} \nonumber \\[-12pt]
    &+ \sum_{s \in S} b_{i,s} g_{i,s,t,a} \Big)
    \label{gen_def} \\[-2pt]
    E^{\text{total}} := \sum_{a \in A} \pi_{a} \sum_{i \in I} \sum_{t \in T} 
    \Big( E_{i}^{\text{min}} u_{i,t,a} &+ {E}_{i}^{\text{su}} v_{i,t,a} \nonumber \\[-12pt]
    &+ \sum_{s \in S} h_{i,s} g_{i,s,t,a} \Big)
    \label{emissions_def}
\end{align}

\vspace{-8pt}

subject to:

\vspace{-18pt}
\begin{gather}
    g_{i,t,a} = g_{i}^{\text{min}} u_{i,t,a} + \sum_{s \in S} g_{i,s,t,a} \hspace{2pt} ; \hspace{2pt} \forall i \in I, t \in T, a \in A  \label{gen_sum} \displaybreak[0] \\[-2pt]
    0 \le g_{i,s,t,a} \le g_{i,s}^{\text{max}} u_{i,t,a} \hspace{8pt} \forall i \in I, s \in S, t \in T, a \in A
    \label{seg_limit} \displaybreak[0] \\
    v_{i,t,a} + z_{i,t,a} \le 1 \hspace{2pt} ; \hspace{2pt} \forall i \in I, t \in T, a \in A
    \label{binary_1} \displaybreak[0] \\
    v_{i,t,a} - z_{i,t,a} = u_{i,t,a} - u_{i,t-1,a} \hspace{2pt} ; \hspace{2pt} \forall i \in I, t \in T, a \in A
    \label{binary_2} \displaybreak[0] \\
    \sum_{\tau=t-\text{g}_{i}^{\text{up}}+1}^{t} v_{i,\tau,a} \le u_{i,t,a} \hspace{2pt} ; \hspace{2pt} \forall t \in T, i \in I, a \in A
    \label{updown_1} \displaybreak[0] \\[-2pt]
    \sum_{\tau=t-\text{g}_{i}^{\text{down}}+1}^{t} z_{i,\tau,a} \le 1 - u_{i,t,a} \hspace{2pt} ; \hspace{2pt} \forall t \in T, i \in I, a \in A
    \label{updown_2} \displaybreak[0] \\
    - r_{i}^{\text{down}} \le g_{i,t,a} - g_{i,t-1,a} \le r_{i}^{\text{up}}
    \hspace{2pt} ; \hspace{2pt} \forall t \in T, i \in I, a \in A \label{ramp}
\end{gather}

\vspace{-22pt}

\begin{multline}
    \sum_{i \in I} m_{i,b}^{\text{unit}} g_{i,t,a} - \sum_{l \in L} m_{l,b}^{\text{line}} f_{l,t,a} - s_{b,t,a}^{\text{ren}} = \\[-3pt]
    d_{b,t,a} - {s}_{b,t,a}^{\text{load}}
    \hspace{2pt} ; \hspace{2pt} \forall b \in B, t \in T, a \in A \label{balance}
\end{multline}

\vspace{-22pt}

\begin{gather}
    -f_{l}^{\text{max}} \le f_{l,t,a} \le f_{l}^{\text{max}} \hspace{2pt} ; \hspace{2pt} \forall l \in L, t \in T, a \in A \label{flow_lim} \displaybreak[0] \\
    f_{l,t,a} = \frac{1}{x_{l}} \sum_{b \in B} m_{l,b}^{\text{line}} \theta_{b,t,a},  \hspace{2pt} ; \hspace{2pt} \forall l \in L, t \in T, a \in A \label{flow_def}
\end{gather}

\vspace{-12pt}

\begin{align}
    \sum_{i \in I \setminus R} u_{i,t,a} (g_{i}^{\text{max}} - g_{i,t,a}) \ge &3\% \sum_{b \in B} d_{b,t,a} + 5\% \sum_{i \in R} g_{i,t,a} \nonumber \\[-2pt]
    + &\max_{i \in I} g_{i}^{\text{max}} \hspace{2pt} ; \hspace{4pt} \forall t \in T, a \in A \label{reserve_capacity}
\end{align}

\vspace{-16pt}

\begin{align}
    \sum_{i \in I} \min \left( r_{i}^{\text{up}} u_{i,t,a}, \left( g_{i}^{\text{max}} - g_{i,t,a} \right) \right) \ge w_{t,a}^{\text{up}} + d_{t,a}^{\text{ramp}} \hspace{2pt}; \label{ramp_suf_1} \nonumber \\[-10pt]
    \forall t \in T, a \in A \displaybreak[0] \\
    \sum_{i \in I} \min \left( \left( r_{i}^{\text{down}} u_{i,t,a}, g_{i,t,a} - g_{i}^{\text{min}} \right) \right) \ge w_{t,a}^{\text{down}} + d_{t,a}^{\text{ramp}} \hspace{2pt};  \label{ramp_suf_2} \nonumber \\[-10pt]
    \forall t \in T, a \in A \displaybreak[0] \\
    0 \le s_{b,t,a}^{\text{load}} \le d_{b,t,a} \hspace{2pt} ; \hspace{4pt} \forall b \in B, t \in T, a \in A \label{shed_limit} \displaybreak[0] \\
    0 \le s_{b,t,a}^{\text{ren}} \le \sum_{i \in R} m_{i,b}^{\text{unit}} g_{i,t,a} \hspace{2pt} ; \hspace{4pt} \forall b \in B, t \in T, a \in A \Big\} \label{spill_limit}
\end{align}
\vspace{-4pt}

The regulator's objective is given in \eqref{upper_obj}, and constrained by the emission limit \eqref{te_max} and the lower-level problem \eqref{lower_obj}-\eqref{spill_limit}.  The system operator's objective is given in \eqref{lower_obj}-\eqref{emissions_def}.  Generator costs curves are piecewise linear \eqref{gen_sum}-\eqref{seg_limit}.  Binary commitment variables are defined in \eqref{binary_1}-\eqref{binary_2} and generator minimum up- and down-times are constrained using \eqref{updown_1}-\eqref{updown_2}.  Generator ramp rate constraints are given in \eqref{ramp}. Power balance is given by \eqref{balance}.  Line flow limits are given by \eqref{flow_lim}-\eqref{flow_def}.  Operating reserve requirements based on the 3+5\% and $N\text{-1}$ policies are ensured using \eqref{reserve_capacity} with flexibility requirements ensured in \eqref{ramp_suf_1}-\eqref{ramp_suf_2}.  Load and renewable generation shedding is constrained by physical limits in \eqref{shed_limit}-\eqref{spill_limit}.

For the intertemporal constraints \eqref{binary_2}-\eqref{ramp}, time periods before the first are treated cyclically.  For instance, $t=24$ is substituted for $t=0$, and $t=23$ for $t=-1$.  This ensures that end-of-day commitments are feasible and that initial conditions are representative, assuming that the days surrounding the representative day are substantially similar. Considering constraints \eqref{reserve_capacity}-\eqref{ramp_suf_2}, the ideal quantity of regulation and load-following reserves is an active research topic \cite{ela2011,fernandezblanco2017}; for simplicity, we use the heuristic 3+5\% rule originally proposed in \cite{wwsis2010}. This formulation assumes a perfectly competitive market; otherwise, the impact of the carbon tax on emissions may vary, as shown in \cite{dearce2016}.

\section{Solution Techniques} \label{taxrate}

One approach to finding the optimal tax rate would be to solve a standard unit commitment over the set of representative days, with a constraint on total emissions, and to take the marginal value of the emissions constraint as the tax rate:

\vspace{-10pt}
\begin{gather}
    \min C^{\text{gen}} + C^{\text{shed}} \label{simple_obj} \\
    \text{Equations \eqref{shed_def}-\eqref{spill_limit}} \label{generic_constraints2} \displaybreak[0] \\
    E^{\text{total}} \le E^{\text{max}} : \lambda \label{TE_constraint}
\end{gather}

\vspace{-4pt}
\noindent
where $\lambda$ denotes the marginal value of the constraint.

\begin{figure}[h]
    \centering
    \includegraphics[width=0.6\linewidth]{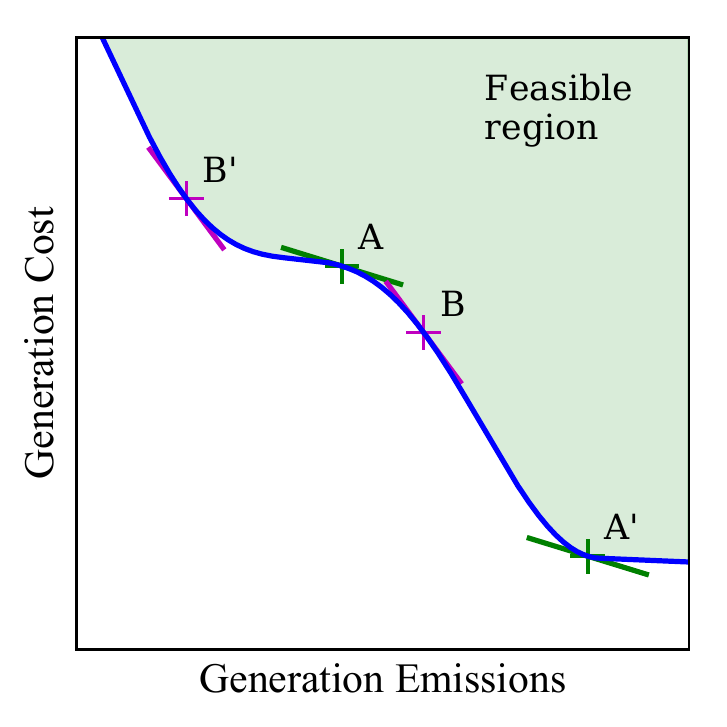}
    \caption{The marginal values at points \{A,B\}, when used as tax rates, result in the solutions at \{A',B'\}.}
    \label{nonconvex_pareto}
\end{figure}

We call such an approach ``Constrained Emission Marginal Value (CEMV)'' method.  However, due to the non-convexity of the UCCT problem (due to binary variables), this approach is liable to produce sub-optimal solutions.  Varying $E^{\text{max}}$ will find solutions on the Pareto frontier of the feasible cost/emissions space, but the resulting $\lambda$, when used as $P^{\text{CO}_2}$ in the UCCT, may find different solutions.  This is because concave portions of the Pareto frontier may not be found by the linearly weighted UCCT formulation, since the optima only exist on the convex hull of the Pareto frontier \cite{deb2016}.  This undesirable outcome is illustrated in Fig. \ref{nonconvex_pareto} and demonstrated for the test system in Section \ref{results}.  Depending on where in the convex region $E^{\text{max}}$ falls, the value of $\lambda$ (\textit{i.e.} the slope of the curve) when input as $P^{\text{CO}_2}$ into the UCCT problem may find: a) a solution with emissions which are greater than the target (A $\rightarrow$ A' in the figure), b) a solution in which emissions are lower than the target, but the production cost is higher than necessary (B $\rightarrow$ B'), c) the optimal cost/emissions point, but at a $P^{\text{CO}_2}$ that is larger than necessary, or d) the minimum $P^{\text{CO}_2}$ which results in the optimal cost/emissions point.  Therefore, it should not be assumed that the CEMV method can find (d) reliably.

\begin{figure}[h]
    \centering
    \includegraphics[width=\linewidth]{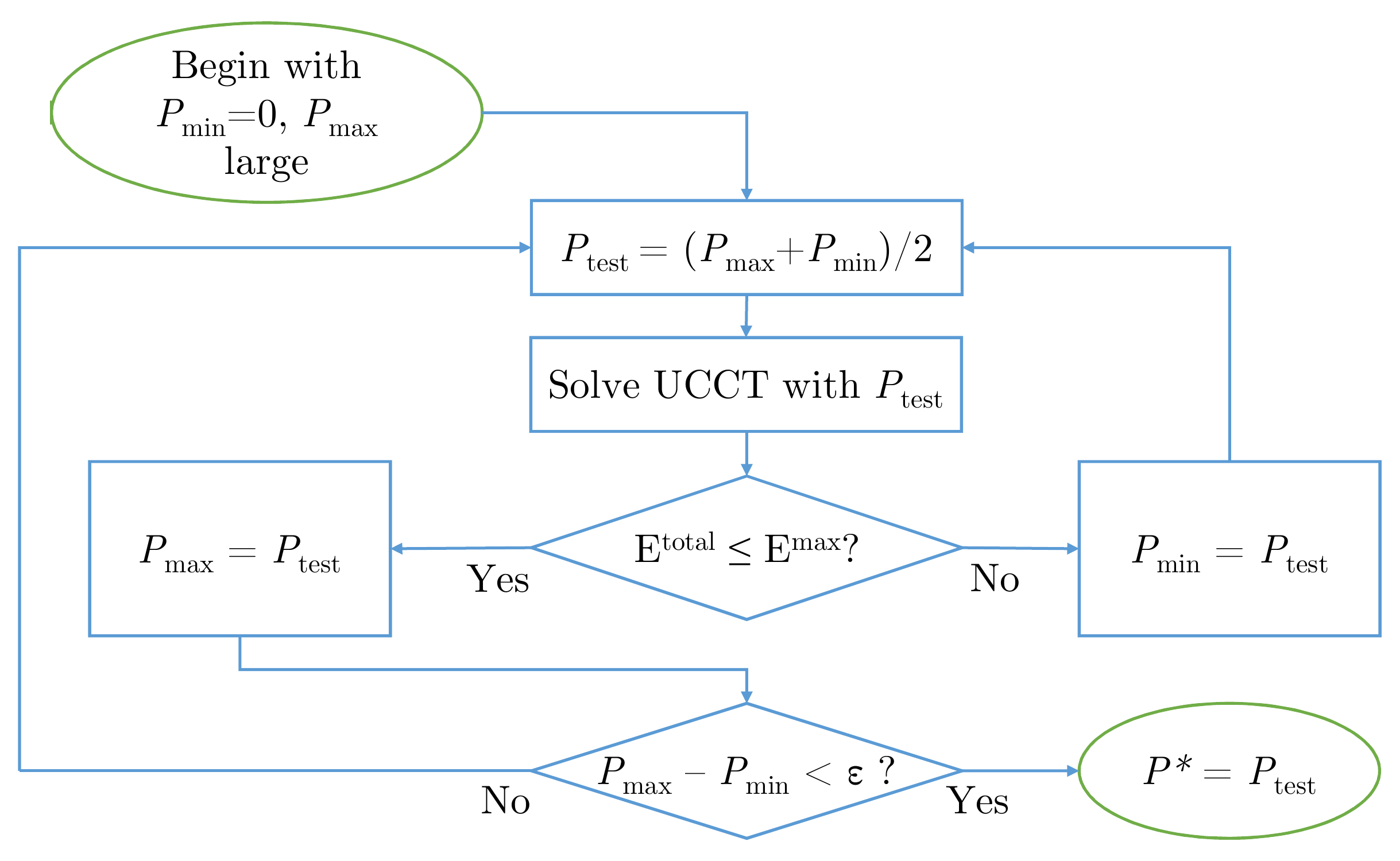}
    \caption{Flowchart for finding optimal $P^{\text{CO}_2}$ using the Weighted Sum Bisection method.}
    \label{flowchart}
\end{figure}

By contrast, the WSB method finds the optimal tax rate by iteratively guessing a $P^{\text{CO}_2}$ value, solving the UCCT problem, and tuning $P^{\text{CO}_2}$ using the bisection method.  Briefly, if a zero of a continuous function is known to be in a certain interval, it can be reliably found by repeatedly bisecting the interval and selecting the sub-interval in which the root must lie, based on the sign of the function value at the midpoint.  Since the goal is to find the tax rate resulting in emissions at or below a certain target, the function is $f(P^{\text{CO}_2}) \hspace{-1pt} = \hspace{-1pt} E^{\text{total}}(P^{\text{CO}_2}) \hspace{-1pt} - \hspace{-1pt} E^{\text{max}}$, and its zero-crossing is at the value of the optimal tax rate, where $E^{\text{total}}(P^{\text{CO}_2})$ is found for a given value of $P^{\text{CO}_2}$ by solving the UCCT.  For a given value of $P^{\text{CO}_2}$, the UCCT can be solved independently for each representative day, aiding computation.  For a non-convex Pareto frontier such as ours, the values of the individual objectives as a function of the weighting factor are noncontinuous but monotonic. Therefore, the WSB method is guaranteed to find the smallest tax rate resulting in emissions at or below the target, if this target is feasible. A very high tax rate (e.g. \$1,000/ton) can be used to estimate the maximum feasible emissions reduction and set the upper bound of the tax range. Since precision is doubled with each iteration, convergence is linear \cite{wood1992}. This approach is shown in Fig. \ref{flowchart}.

\section{Case Study} \label{casestudy}

The electrical system for this case study is a modified  ISO New England (NE) test system \cite{krishnamurthy2016}.  Data from the Energy Information Administration (EIA) are used for fuel prices \cite{eia_prices} and for per-MMBTu $\text{CO}_{2}\text{e}$ emissions by fuel \cite{eia_co2}.  Though variability from renewable generation can induce additional $\text{CO}_2$ emissions from thermal generators \cite{katzenstein2009}, this effect was not modeled in the case study. Five representative days are chosen using a hierarchical clustering algorithm \cite{pitt2000} and run at a one-hour time resolution.  The load shed penalty is set at \$10,000/MWh and the renewable spillage penalty is \$20/MWh.  Ramping requirements are set such that the system has the capacity to react to 1\%/hour load ramps and all wind farms ramping their production $\pm$ 20\% over one hour, based on analysis of Bonneville Power Administration wind power production data in \cite{nosair2015}.  This case study was implemented using GAMS v24.0 and solved using CPLEX v12.5 with a 0.1\% optimality gap on an Intel Xenon 2.55 GHz processor with at least 32 GB RAM.

\section{Results} \label{results}

\subsection{System Characteristics}

\begin{figure}[h]
    \centering
    \includegraphics[width=\linewidth]{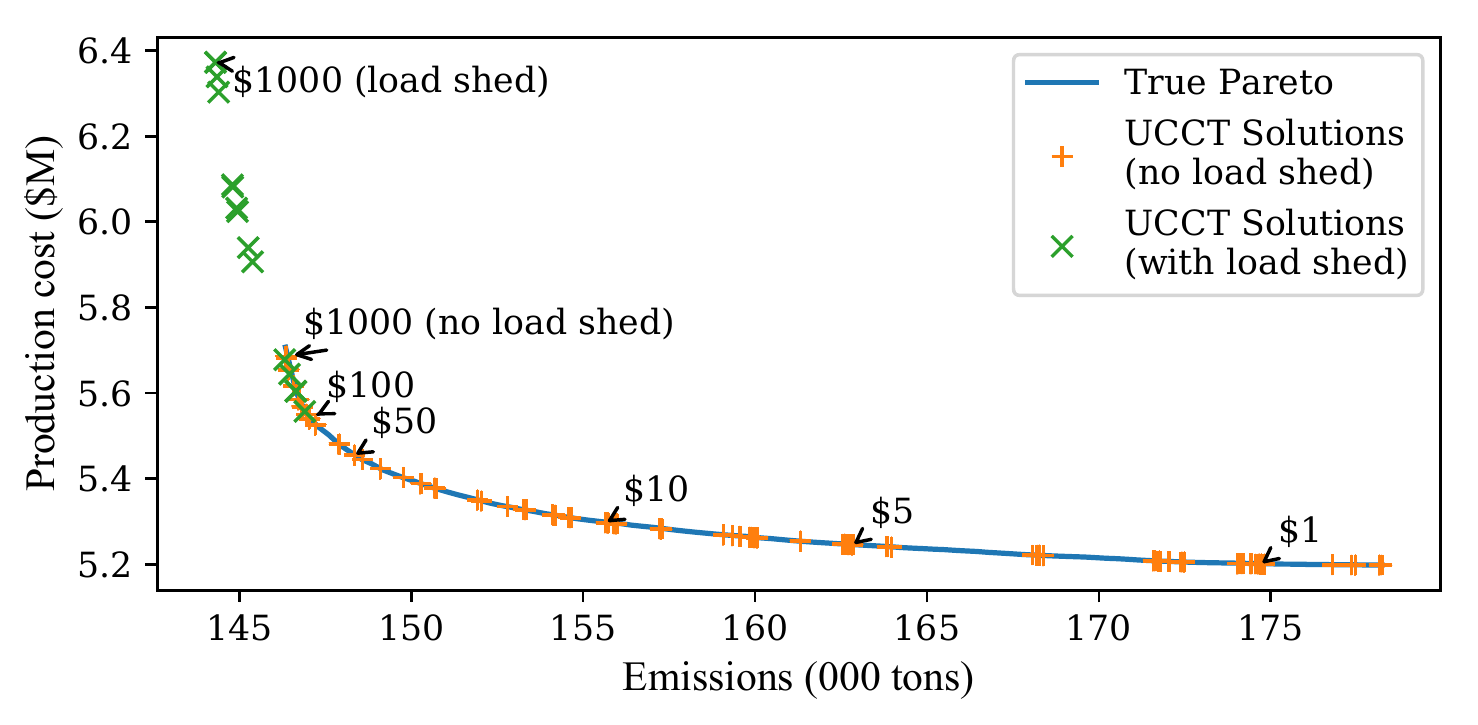}
    \caption{The cost/emissions Pareto frontier. The line is points found by constraining emissions, the crosses are points found by varying $P^{\text{CO}_2}$.}
    \label{three_paretos}
\end{figure}

\begin{figure}[h]
    \centering
    \includegraphics[width=\linewidth]{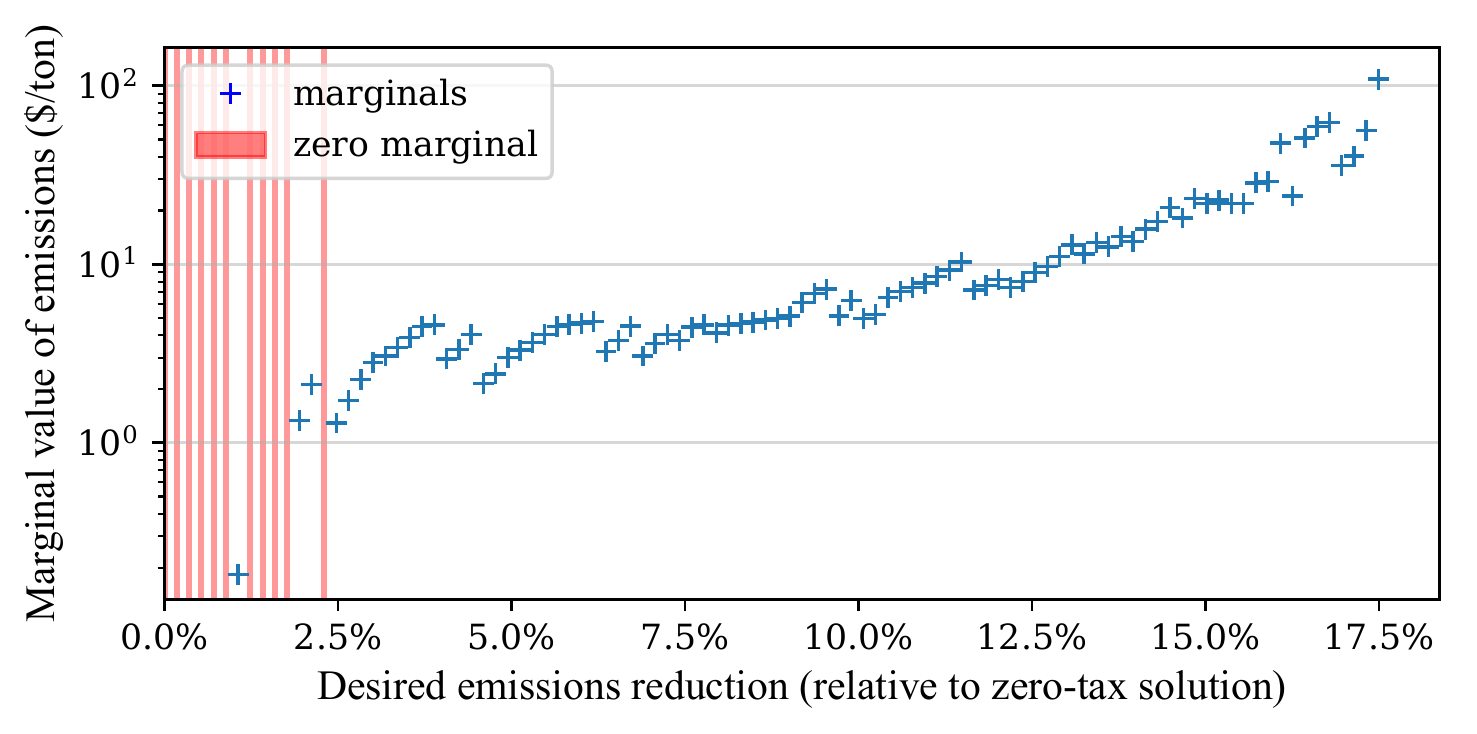}
    \caption{Marginal value found at 100 sample points on the Pareto frontier.}
    \label{discovered_prices}
\end{figure}

Fig. \ref{three_paretos} shows the Pareto frontier of the trade-off between emissions and production costs (\textit{i.e.} fuel and shed costs).  The full Pareto frontier is sampled at 100 equally-spaced points by constraining emissions and varying $E^{\text{max}}$, and the convex hull of the cost/emissions space is sampled by using the UCCT and varying $P^{\text{CO}_2}$, with and without load-shedding.  Load shedding is only economically justified under very high tax rates and results in very high costs, so load-shedding solutions are omitted in all following figures for the sake of clarity.  Though the Pareto frontier may at first glance appear convex, there are many small concave regions.  This can be seen by plotting the marginal value at the sample points, as shown in Fig. \ref{discovered_prices}; since the marginal values do not increase monotonically, the frontier must be non-convex \cite{boyd2004}.

\subsection{Determining a Tax Rate to Meet Policy Goals}

\begin{figure}[h]
    \centering
    \includegraphics[width=\linewidth]{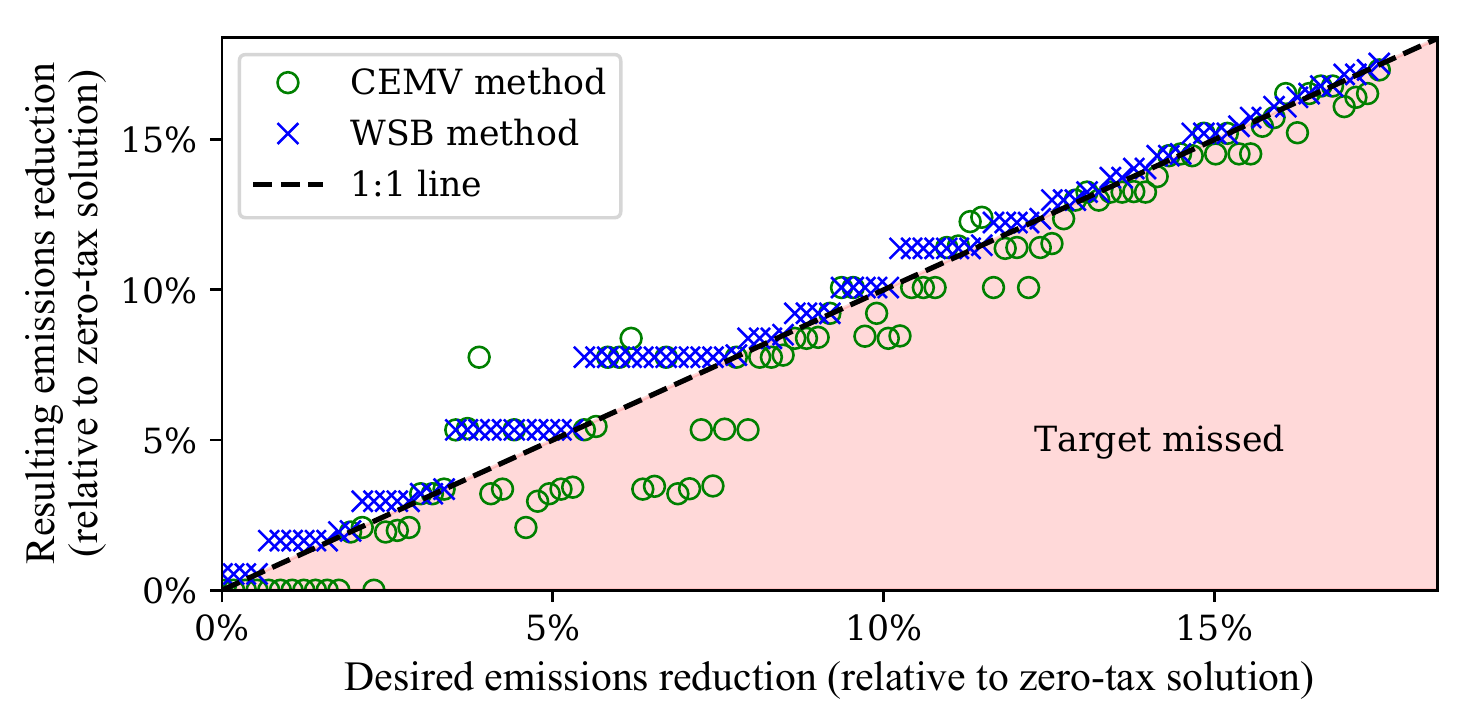}
    \caption{Comparison of results using the CEMV method and the WSB method.}
    \label{marginal_v_heuristic}
\end{figure}

\begin{figure}[h]
    \centering
    \includegraphics[width=\linewidth]{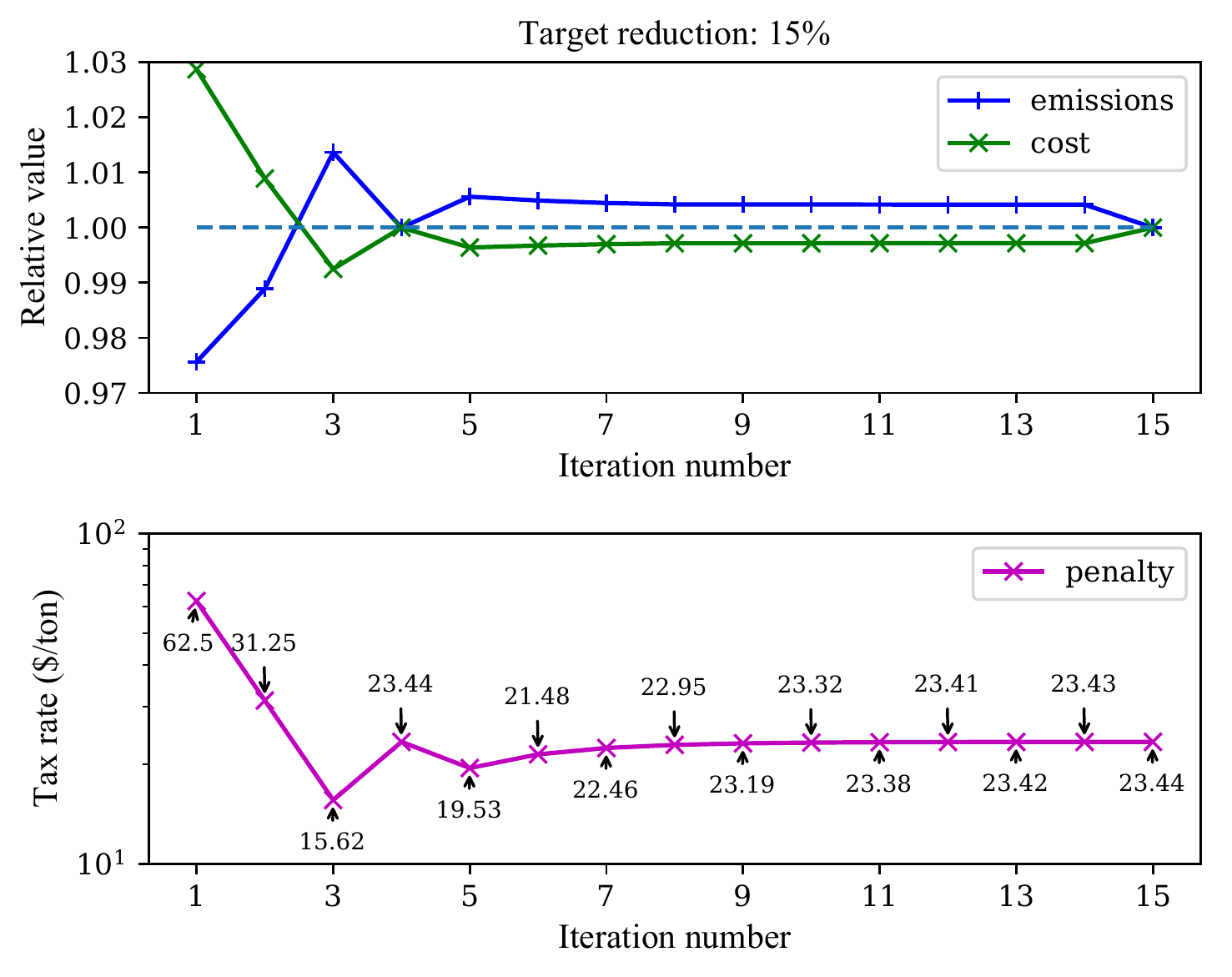}
    \caption{Convergence of WSB method to final value.}
    \label{bisection_demo}
\end{figure}

If the CEMV method were used to set a tax rate, there is no guarantee that the solution to the UCCT problem would meet the desired emissions reduction.  This is illustrated in Fig. \ref{marginal_v_heuristic}, which uses the same set of sample points as Figs. \ref{three_paretos} and \ref{discovered_prices}.  As shown, values of $P^{\text{CO}_2}$ derived from the CEMV method do not reliably meet their desired emissions reductions when used in the UCCT.  By comparison, the WSB method is guaranteed to meet or exceed the emissions reduction target.  Convergence of the WSB method to its final values is shown in Fig. \ref{bisection_demo} for a target emissions reduction of 15\%.  The WSB method, given an emissions target, reliably converges to an optimal tax rate within $1\cent$ from an initial range of \$0-\$100/ton in 14 iterations of the UCCT problem.  Though the tax rate which is converged upon may not be the true optimum, it can be shown that the solution exceeds the true optimum by no more than a specified tolerance, and this tolerance can be halved with each additional iteration of the UCCT problem. By comparison, naively finding the tax rate by solving for each possible rate in $1\cent$ increments would require 10,000 solves. The wider the range of potential solutions, and the greater the desired accuracy, the more efficiently the WSB performs.

\begin{figure}[h]
    \centering
    \includegraphics[width=\linewidth]{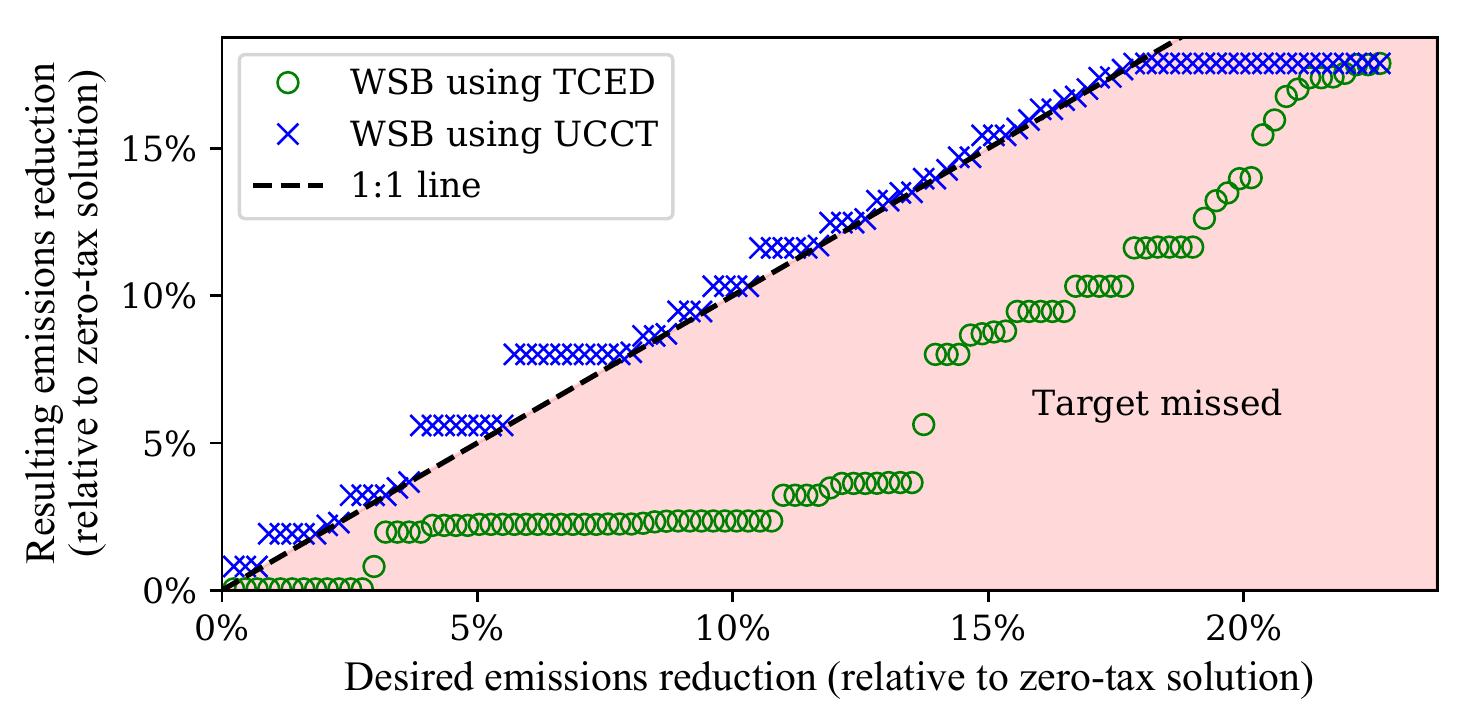}
    \caption{Comparison of WSB results when ignoring binary variables and intertemportal constraints (TCED problem) vs. including them (UCCT method).}
    \label{TCED_vs_WSB}
\end{figure}

The importance of including binary variables is illustrated in Fig. \ref{TCED_vs_WSB}.  For this figure, the UCCT formulation is transformed into a transmission-constrained economic dispatch (TCED) problem by ignoring intertemporal constraints \eqref{binary_2}-\eqref{ramp} and setting $g_{i}^{\text{min}},C_{i}^{\text{min}},E_{i}^{\text{min}}=0$ for all generators.  By using the WSB method to find the required tax rate for a given desired emissions reduction for the TCED problem, and inputting that resulting tax rate into the UCCT problem, it can be seen that the realized emissions reductions fail to meet the targets. Factors which can contribute to this outcome include: the requirement to burn fuel to synchronize generators on start-up, and the requirement to commit additional generators to prepare for large ramps in net load, which occur more commonly and with greater magnitude with the introduction of large quantities of renewables.

\subsection{Handling Uncertainty}

Since this formulation requires estimating the distribution of representative days in a future year, there is some uncertainty in the actual realization.  The variance in annual realized emissions, based on this sampling probability, is given by \eqref{weather_uncertainty}.  If policy-makers desire to achieve emissions reductions with a specified level of certainty, tax rates can be set such that the likelihood of achieving such reduction happens with the desired probability using \eqref{cdf} due to the Central Limit Theorem \cite{billingsley1995}.

\vspace{-12pt}
\begin{gather}
    \text{Var}[E^{\text{total}}] = \sigma_{E}^2 = 365 \sum_{a \in A} \pi_{a} (E_{a} - E^{\text{total}})^2
    \label{weather_uncertainty} \displaybreak[0]\\[-2pt]
    \text{Prob}[E^{\text{total}} \le E^{\text{max}}] = \Phi \left( \frac{E^{\text{max}}-E^{\text{total}}}{\sigma_{E^{\text{total}}}}\right) \label{cdf}
\end{gather}
\vspace{-12pt}

\noindent
where $\Phi(\cdot)$ is the cumulative distribution function of the standard normal distribution.

\begin{figure}
    \centering
    \includegraphics[width=\linewidth]{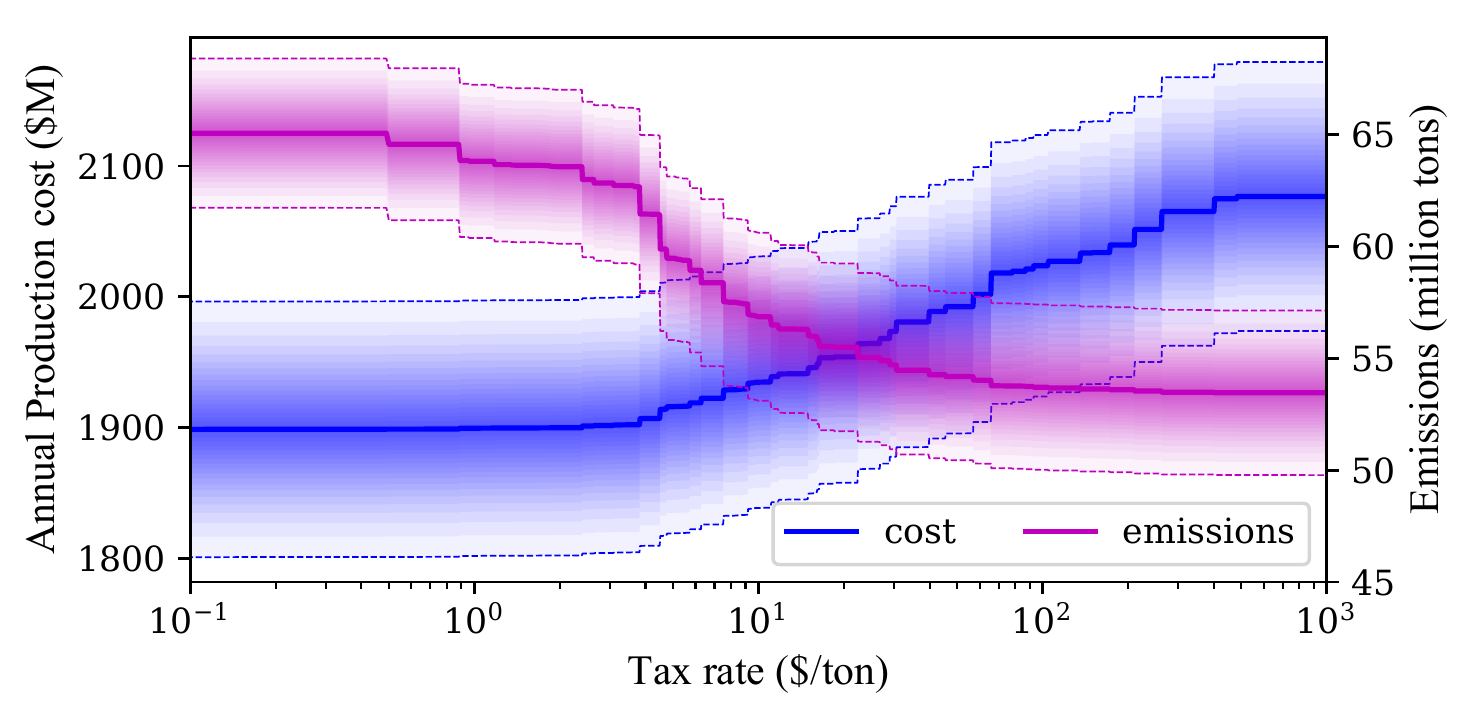}
    \caption{Fuel cost and emissions as a function of tax rate, incorporating weather uncertainty. Bands represent 95\% certainty range.}
    \label{psweep_w_weather_unc}
\end{figure}

\begin{figure}
    \centering
    \includegraphics[width=\linewidth]{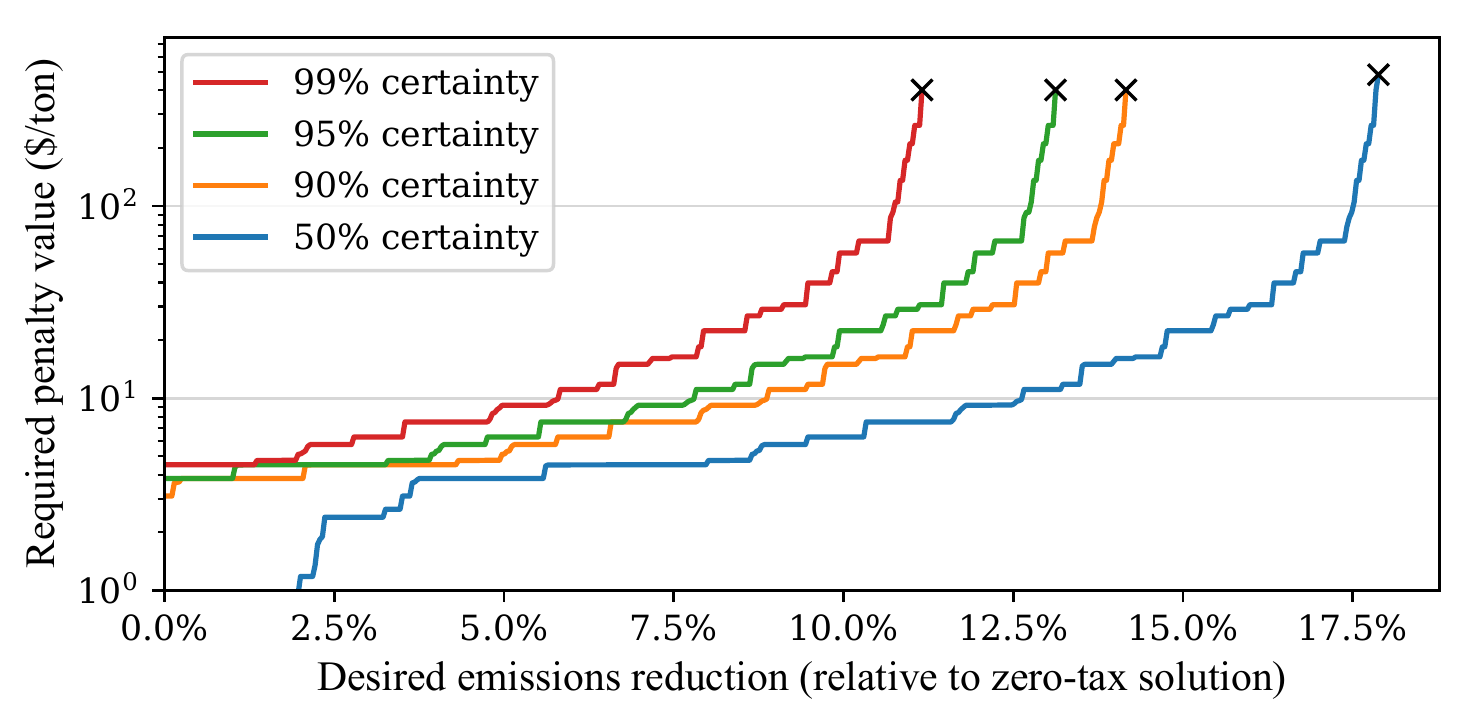}
    \caption{Tax rate required to achieve emissions reductions, based on weather uncertainty.}
    \label{psweep_inverse}
\end{figure}

The choice of increasing or decreasing the tax rate in the WSB method is then based on whether this likelihood meets the desired level of certainty.  Fig. \ref{psweep_w_weather_unc} illustrates the uncertainty range around the expected cost and emissions, and Fig. \ref{psweep_inverse} illustrates the tax rate required to achieve a desired emissions reduction for various values of certainty.  As shown, the required tax rate to meet a given emissions target increases with the level of certainty required, and some emissions reductions targets which are achievable on average are not able to be met with much certainty, no matter the tax rate.

A similar process can be used in order to handle the uncertainty of fuel prices.  Currently in the United States, abundant shale gas makes gas-fired power plants more competitive, but these low prices may not persist.  If there is a desire to set a tax rate to be robust to fluctuations in the price of natural gas, the tax-setting process can be run using the highest gas price that can be reasonably expected.

\begin{figure*}
    \centering
    \includegraphics[width=\textwidth]{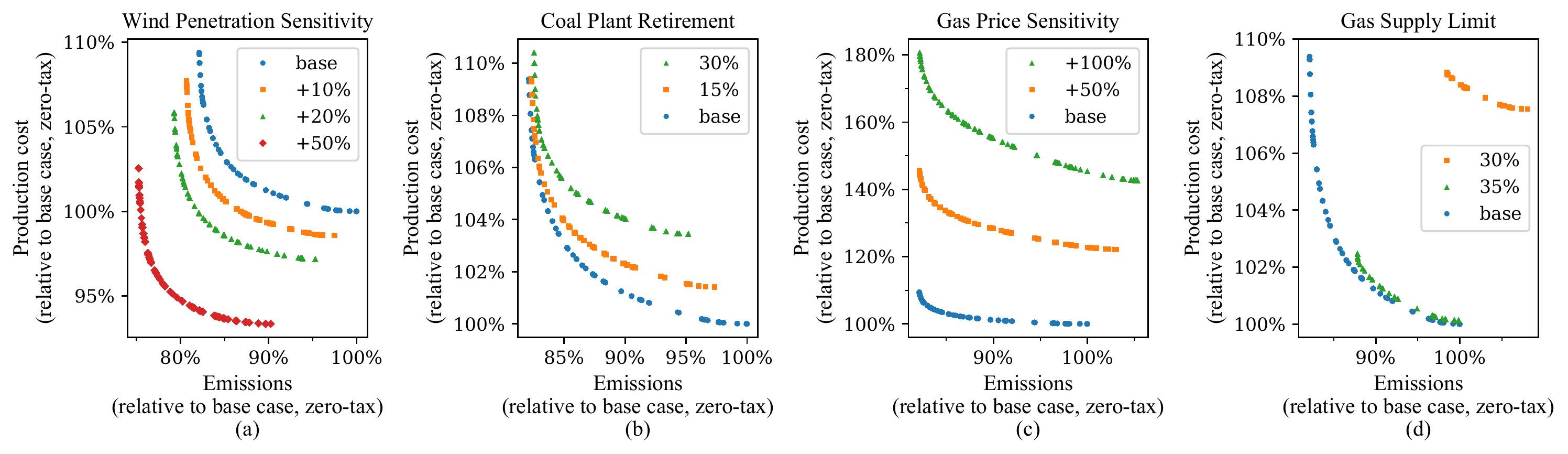}
    \caption{Sensitivities of Pareto frontier to (a) wind penetration, (b) coal plant retirement, (c) natural gas price, and (d) gas supply limit.}
    \label{subplot_pareto_sensitivities_diverse}
\end{figure*}

\begin{figure*}
    \centering
    \includegraphics[width=\textwidth]{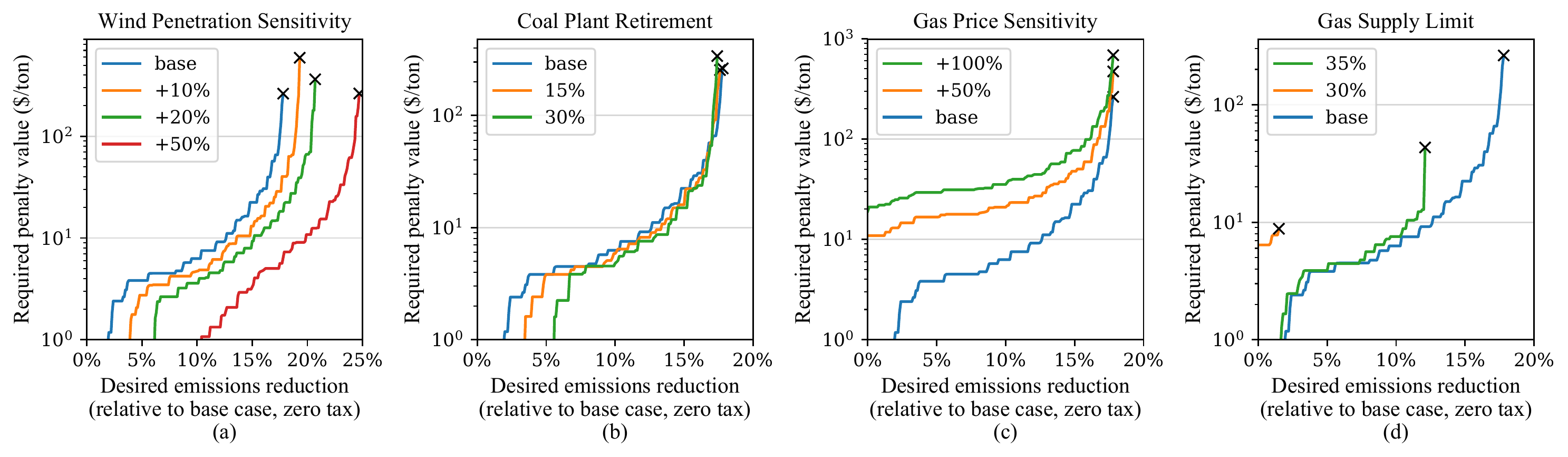}
    \caption{Sensitivities of required tax rate to (a) wind penetration, (b) coal plant retirement, (c) natural gas price, and (d) gas supply limit.}
    \label{subplot_inv_pareto_sensitivities_diverse}
\end{figure*}

\begin{figure*}
    \centering
    \includegraphics[width=\textwidth]{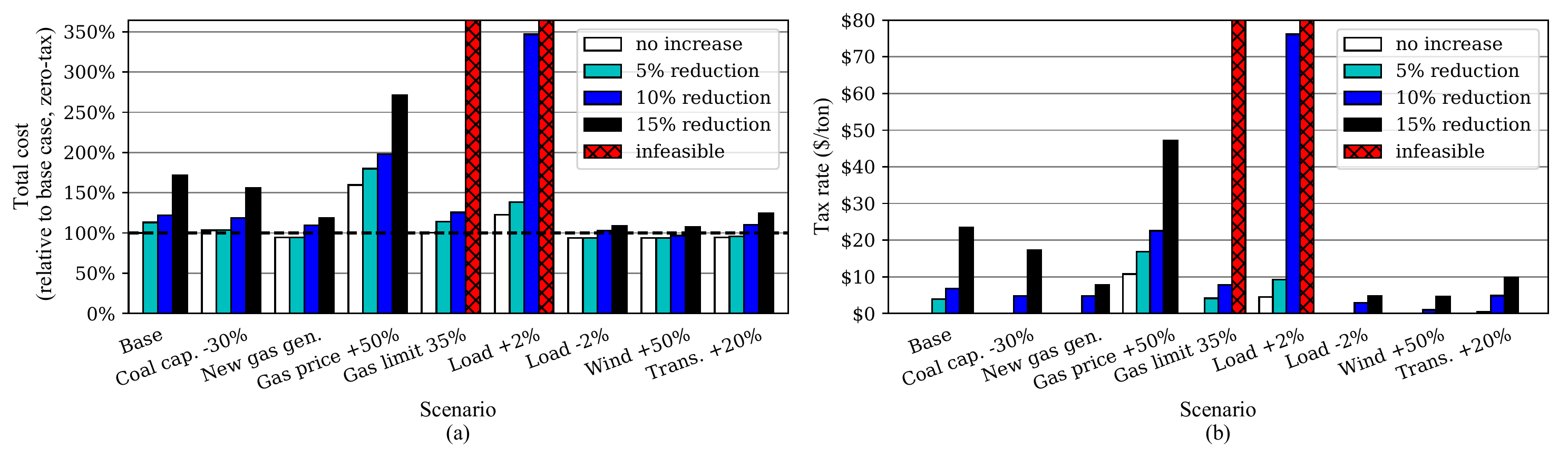}
    \caption{Results of selected scenarios for specified emission reduction targets, (a) total cost, (b) required tax rate.}
    \label{scenario_bar_results}
\end{figure*}

\subsection{Sensitivity Analyses}

Sensitivity analyses to changes in the system's wind penetration, coal plant retirement, gas prices, and gas supply limitations on UCCT solutions are shown in Figs. \ref{subplot_pareto_sensitivities_diverse}(a)-(d), respectively, and the impact on the tax rates required to achieve desired emissions reductions are shown in Figs. \ref{subplot_inv_pareto_sensitivities_diverse}(a)-(d).  Additional scenarios are also run for specific emissions reductions targets and shown in Fig. \ref{scenario_bar_results}:

\begin{itemize}
    \item \textbf{Wind Penetration}: Wind penetration, initially at 8\% of total energy, is increased by 10, 20, or 50\%.
    \item \textbf{Coal Retirement}: Coal generation is retired, either one or two highest-cost generators (15\% or 30\% of the coal-generating capacity), consistent with estimates in \cite{rahmani2016}.
    \item \textbf{New Gas Generator}: One new 250 MW gas generator is added at the bus with highest average locational marginal price (LMP), bus 8, increasing the gas capacity by 2.4\%.
    \item \textbf{Gas Price}: The price of natural gas generation is increased by either 50 or 100\%.
    \item \textbf{Gas Limit}: For each day, gas generators are limited in the amount of energy that they can supply, at either 30 or 35\% of daily total energy. This is intended to simulate regional gas shortages such as those experienced in New England \cite{isone2014} and Southern California \cite{sce2014} in 2014.
    \item \textbf{Load Increase/Decrease}: The demand for electricity at each hour is scaled up or down by 2\%.
    \item \textbf{Transmission Capacity}: The capacity of all transmission corridors is increased by 20\%.
\end{itemize}

Several of these scenarios have similar effects: increases in wind penetration, gas generation capacity, or transmission capacity, or decreases in load.  For all of these scenarios, the zero-tax solutions have lower costs and emissions than the base case, a given emissions target can be met with a lower tax rate, and the maximum emissions reduction is increased.  These effects can be seen in Figs. \ref{subplot_pareto_sensitivities_diverse}(a) and \ref{subplot_inv_pareto_sensitivities_diverse}(a) for the wind penetration case and in Fig \ref{scenario_bar_results} for the other cases.

Other scenarios have differing effects.  For the coal plant retirement case, the zero-tax solution has lower emissions than the base case but is more expensive, since the retired coal generation is replaced by more-expensive gas generation.  However, despite the higher zero-tax cost, a desired emissions reduction can be met with a lower tax rate, as shown in Figs. \ref{subplot_pareto_sensitivities_diverse}(b) and \ref{subplot_inv_pareto_sensitivities_diverse}(b), and some targets can be met at a lower total cost. For example, a 15\% GHG reduction target in the base case requires a tax rate of \$23/ton and a total cost of increase of 71\%, but the same target can be met with a \$17/ton tax rate and 56\% cost increase in the 30\% coal retirement case, as shown in Fig. \ref{scenario_bar_results}.  For the gas price increase cases, the zero-tax solution is both more expensive and higher-emitting than the base case, and a desired emissions reduction requires a higher tax rate, but the range of possible emissions reductions is not affected, as shown in Figs. \ref{subplot_pareto_sensitivities_diverse}(c) and \ref{subplot_inv_pareto_sensitivities_diverse}(c).  For the case where there are limitations on the energy supplied by gas generators, the effect is highly dependent on the limit values and emissions targets, as shown in Figs. \ref{subplot_pareto_sensitivities_diverse}(d) and \ref{subplot_inv_pareto_sensitivities_diverse}(d).  When gas generators are limited to providing no more than 35\% of total energy, there is minimal impact for emissions reductions less than 10\%, but past this point there is limited ability to reduce emissions, and reductions require a higher tax rate.  At a limit of 30\%, the zero-tax solution is 7.6\% more expensive and emits 7.6\% more GHGs when compared to the base case, and only a very modest reduction in emissions (1.5\%) is possible.  This illustrates the impact that gas system constraints can have on emission reductions goals.

\begin{figure}
    \centering
    \includegraphics[width=\linewidth]{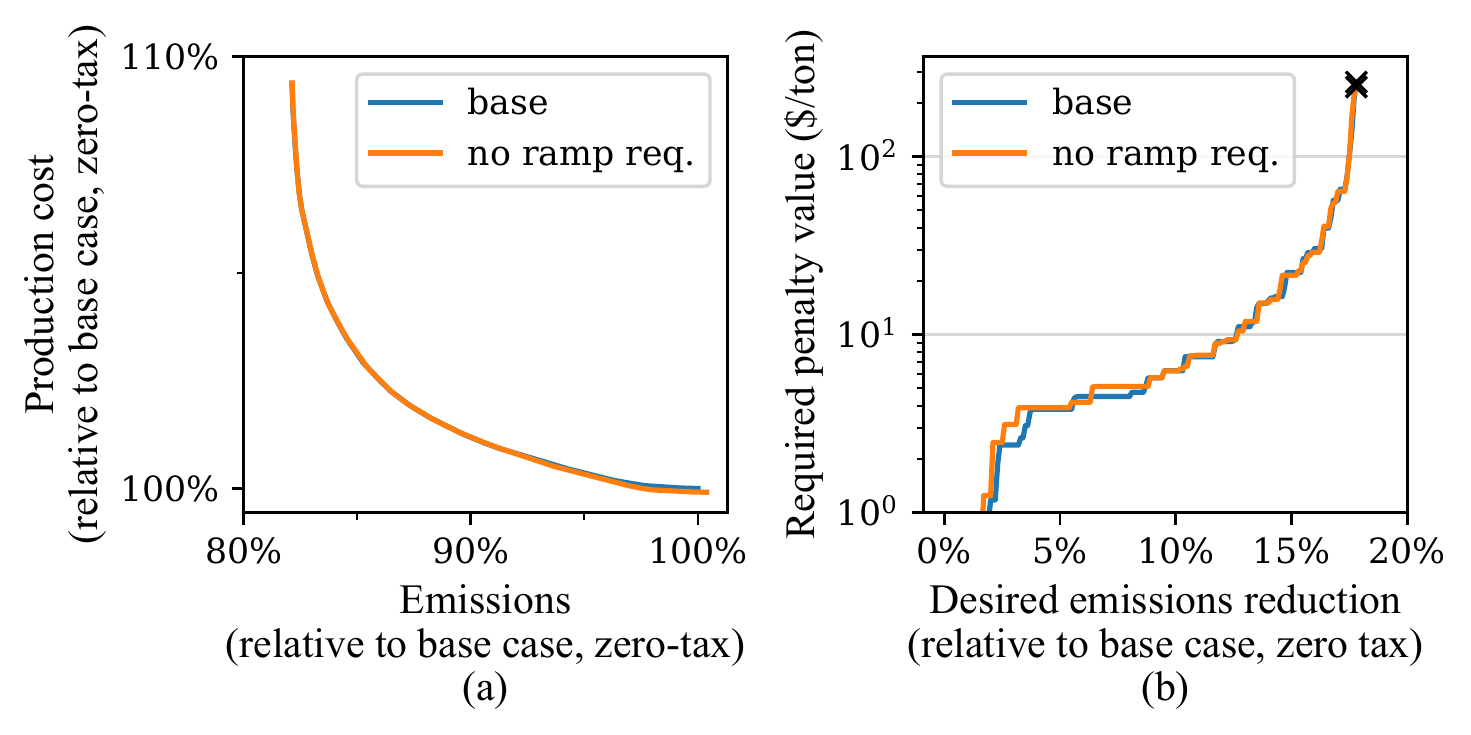}
    \caption{Impact of relaxing system ramp requirement constraint on (a) cost/emissions Pareto frontier, and (b) tax rate required to achieve emissions reductions.}
    \label{ramping_scenario}
\end{figure}

The impact of relaxing the system flexibility constraints \eqref{ramp_suf_1}-\eqref{ramp_suf_2} is also investigated. For this system, the greatest impact is seen at lower tax rates, where relaxing the ramping capability requirement results in slightly higher emissions at slightly lower production cost, as shown in Fig. \ref{ramping_scenario}(a). At zero tax rate, emissions are 0.4\% higher and production cost is 0.1\% lower as compared to the base case. The impact on the tax rate required for a given emission reduction is similar, as shown in Fig. \ref{ramping_scenario}(b): greater impact at low emissions reductions targets, with differences diminishing at more aggressive emissions reductions targets. The converse is true when observing the impact on profit by generation technology: low tax rates have minimal impact, while higher tax rates have more significant impacts. At tax rates of \$0-10/ton, the average impact on profit is within 1.5\% for all generation technologies, while for tax rates of \$10-100/ton, relaxing the flexibility constraints results in profits on average 102\% higher for coal generators, and 5-6\% lower for all other generation technologies.

\subsection{Impact on Investment Decisions}

\begin{figure}
    \centering
    \includegraphics[width=\linewidth]{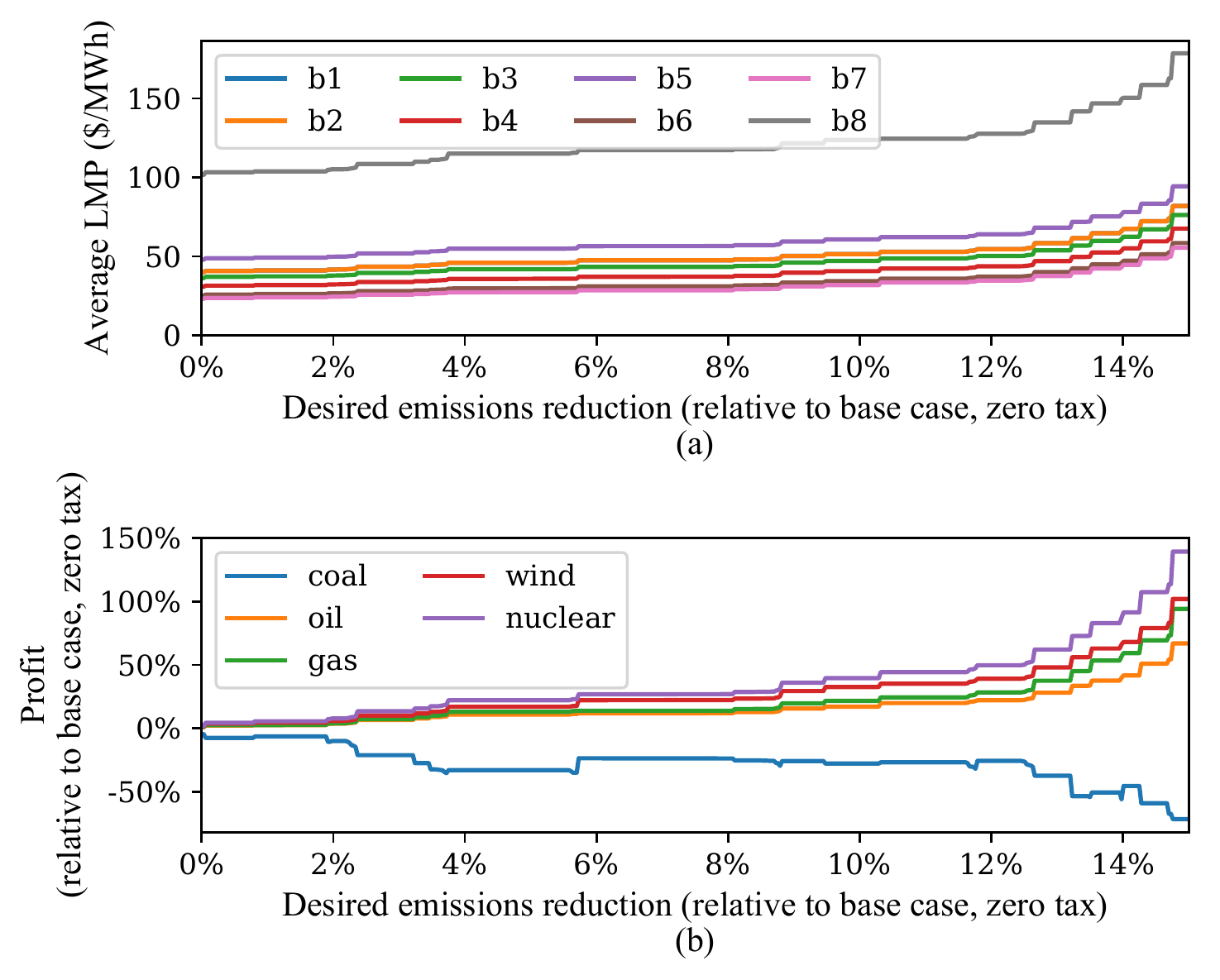}
    \caption{(a) Average LMPs for buses b1-b8 as a function of desired emissions reduction, and (b) Profit by fuel as a function of desired emissions reduction.}
    \label{lmps_and_profit_by_reduction}
\end{figure}

As can be seen in Fig. \ref{lmps_and_profit_by_reduction}(a), larger emissions reduction targets lead to higher average LMPs, which translates to a better value proposition for investing in new non-coal generation (as shown in Fig. \ref{lmps_and_profit_by_reduction}(b), the carbon tax reduces profit by coal generators). For example, at a target emissions reduction of 5\%, the profit for wind generators is 17\% higher as compared to the zero-tax solution; at a target of 10\%, the profit is 33\% higher. In the New Gas Generator scenario, the new generator makes 12\% more profit at an emissions reduction target of 10\%; at a target reduction of 5\%, no tax is required.

Along with higher LMPs, emission reduction targets also result in increased congestion surplus, improving the value proposition for investments in new transmission.  At a target GHG reduction of 5\%, the congestion surplus is increased by 9.5\%; at a 10\% target the surplus is increased by 16.1\%. Congestion surplus as a share of total cost remains relatively constant, however: 20.6\% in the zero-tax case, 20.4\% in the 5\% target case, and 19.6\% in the 10\% target case.

In addition to higher average LMPs, there is also an increase in the variability of LMPs, which improves the value proposition for grid-scale energy storage devices, consistent with findings in \cite{desisternes2016}. At a 5\% emissions reduction target, the average LMP for each bus is 35-60\% higher (average of 46\%), while the standard deviation of LMPs for each bus is 101-113\% higher.  However, while investments in wind and gas generation tend to reduce emissions (and therefore the tax rate required to meet emissions targets), the impact of energy storage is much less clear \cite{lin2016_apen}.

\subsection{Coal Generators and Market Share}

\begin{figure}
    \centering
    \includegraphics[width=\linewidth]{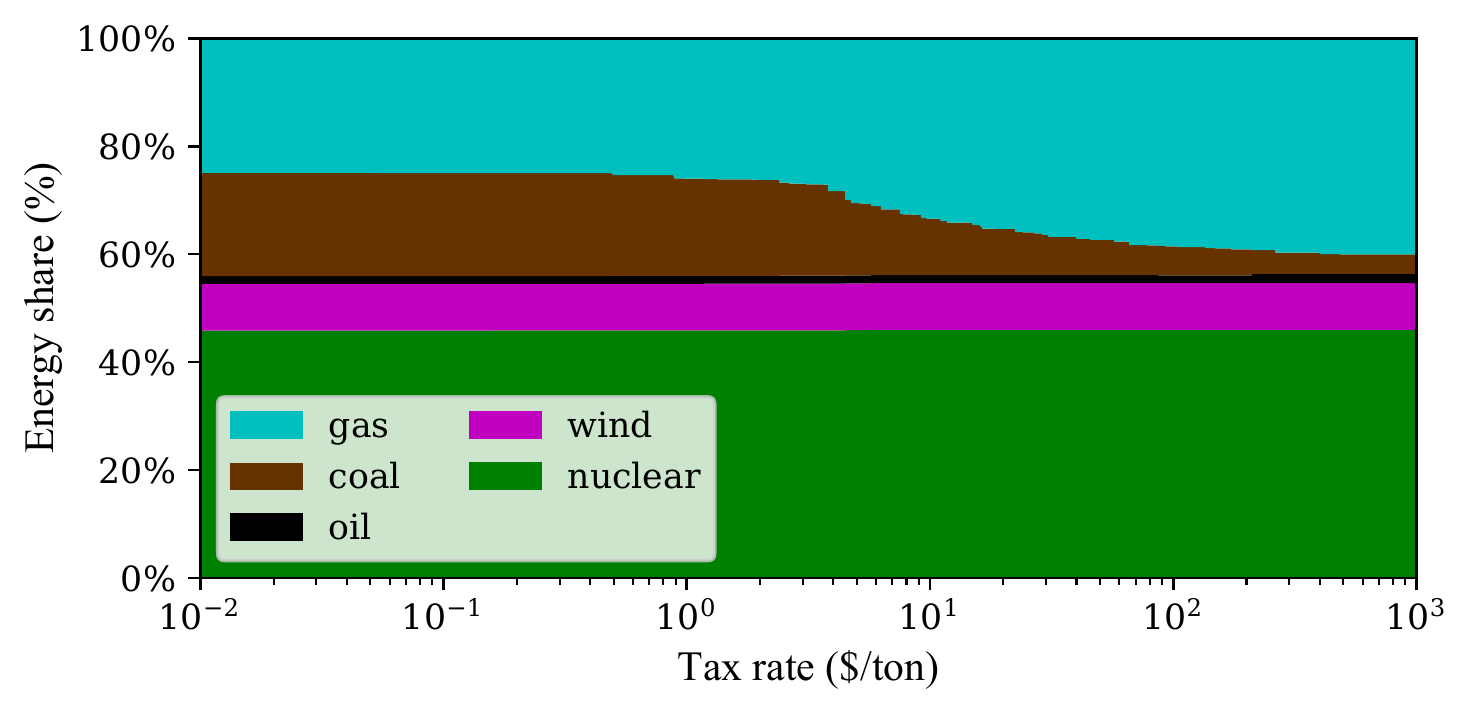}
    \caption{Energy by fuel as a function of tax rate.}
    \label{fuel_areas}
\end{figure}

As would be expected, as the desired emissions reduction and the tax rate required to achieve this reduction both increase, the share of energy which is provided by coal and the profit made by coal generators both decrease, as can be seen in Figs. \ref{lmps_and_profit_by_reduction}(b) and \ref{fuel_areas}. These reduced profits may lead to earlier coal generation retirement and reductions in coal mining employment \cite{heslin1989}.  However, strategic investments in new wind generation and transmission network expansion both allow emissions targets to be met with a lower tax rate compared to the base case, and increase the profit of coal-powered generators, as shown in Table \ref{coal_profits}.   Additionally, emissions taxes which are set based on emissions targets increase the incentive for owners of coal-powered generation to invest in renewable generation, since the increased profits for coal generators caused by a lower tax rate provides an `extra' revenue stream.

\begin{table}
	\caption{Coal profit dependence on investments, 5\% reduction target, given increases in wind penetration or transmission capacity}
    \begin{tabular}{ 
        			>{\centering\arraybackslash}p{0.17\linewidth} >{\centering\arraybackslash}p{0.23\linewidth} >{\centering\arraybackslash}p{0.23\linewidth}
        			>{\centering\arraybackslash}p{0.18\linewidth}
            		}
        \hline \hline
        Scenario & Additional Coal Profit (\$) & Additional Wind Profit (\$) & \% Increase (coal/wind)\\
		\hline +10\% Wind & 24,212 & 77,336 & 31.3\% \\
		\hline +20\% Wind & 62,170 & 140,276 & 44.3\% \\
		\hline +10\% Trans. & 22,835 & N/A & N/A \\
		\hline +20\% Trans. & 131,112 & N/A & N/A \\
		\hline \hline
	\end{tabular}
	\label{coal_profits}
\end{table}

\section{Conclusion} \label{conclusion}

The Weighted Sum Bisection method can be efficiently applied to the problem of determining the optimal tax rate to meet a given emissions target, and has been shown to be compatible with uncertainty in weather, demand, fuel prices, and generation fleet.  In addition to changing short-term generator commitment and dispatch, pricing carbon for emissions targets is also associated with several longer-term effects. 

Higher prices for electricity reduce overall demand and increase the value proposition for investment in new cleaner generation, transmission, and grid-scale energy storage, as well as technologies which generate or consume electricity more efficiently.  Carbon capture and sequestration projects also become more attractive. Conversely, investment in these resources reduces the tax rate required to achieve certain emissions targets.

In addition, carbon pricing collects revenue which can be used to invest in projects to reduce GHG emissions, adapt to climate change, and/or provide economic relief to communities which are negatively impacted by taxation of carbon, such as those heavily reliant on coal-mining.  This revenue stream could also be used to offset tax reductions elsewhere to obtain revenue neutrality.

\bibliographystyle{IEEEtran}
\bibliography{bibliography}

\begin{IEEEbiography}[{\includegraphics[width=1in,height=1.25in,clip,keepaspectratio]{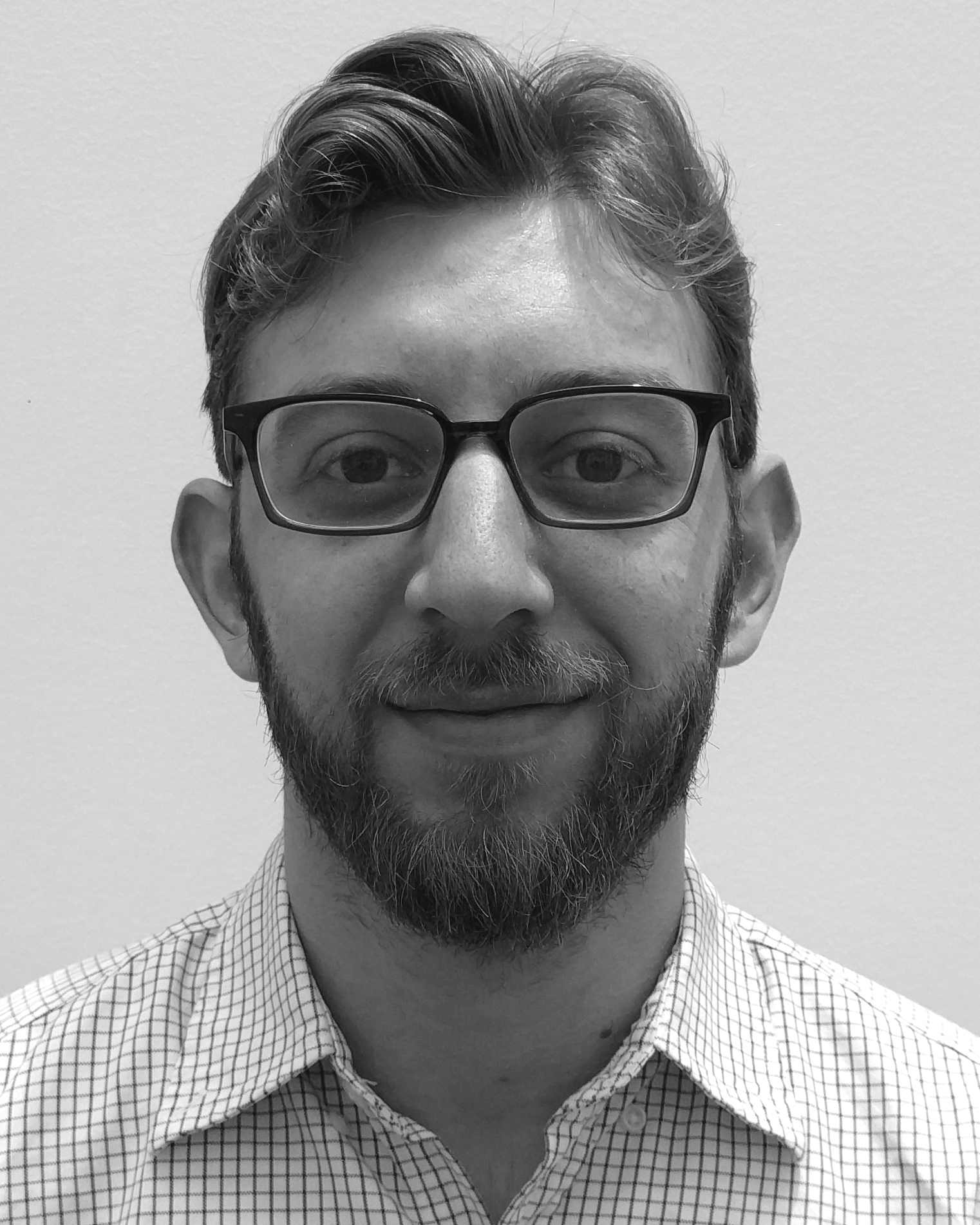}}]{Daniel Olsen}
	(S'14) received the B.Sc. degree in mechanical engineering and electric power engineering from Rensselaer Polytechnic Institute in 2010. He is currently pursuing the Ph.D. degree in electrical engineering at the University of Washington. Previously, he was a Research Associate with Lawrence Berkeley National Laboratory. His research interests include planning and policies for power system emissions, multiple-energy systems, and distributed flexibility resources.
	\end{IEEEbiography}
	
	\begin{IEEEbiographynophoto}{Yury Dvorkin}
	(S’11-M’16) received his Ph.D. degree from the University of Washington, Seattle, WA, USA, in 2016. Dvorkin is currently an Assistant Professor in the Department of Electrical and Computer Engineering at New York University, New York, NY, USA. Dvorkin was awarded the 2016 Scientific Achievement Award by Clean Energy Institute (University of Washington) for his doctoral dissertation “Operations and Planning in Sustainable Power Systems”. His research interests include power system operations, planning, and economics. 
	\end{IEEEbiographynophoto}
	
	\begin{IEEEbiography}[{\includegraphics[width=1in,height=1.25in,clip,keepaspectratio]{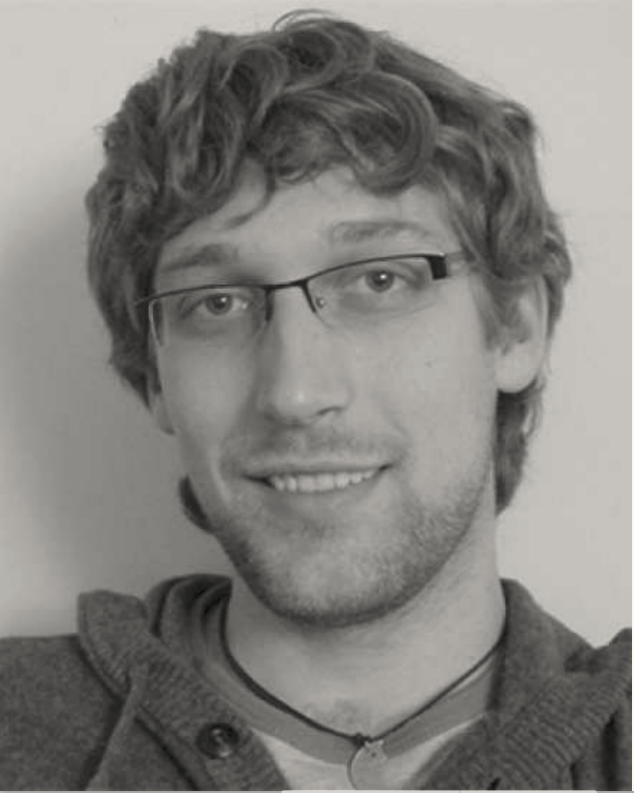}}]{Ricardo Fern\'andez-Blanco}
	(S'10-M'15) received the Ingeniero Industrial degree and the Ph.D. degree in electrical engineering from the Universidad de Castilla-La Mancha, Ciudad Real, Spain, in 2009 and 2014, respectively.
    He is currently a Postdoctoral Researcher at the University of M\'alaga, Spain. He also was a Postdoctoral Researcher at the University of Washington, Seattle, WA, USA. He was working as a Scientific/Technical Project Officer in the Knowledge for the Energy Union Unit at the DG JRC (Joint Research Center) of the European Commission, Petten, The Netherlands. His research interests include the fields of operations and economics of power systems, smart grids, bilevel programming, hydrothermal coordination, electricity markets, and the water-energy nexus.
	\end{IEEEbiography}
	
	\begin{IEEEbiography}[{\includegraphics[width=1in,height=1.25in,clip,keepaspectratio]{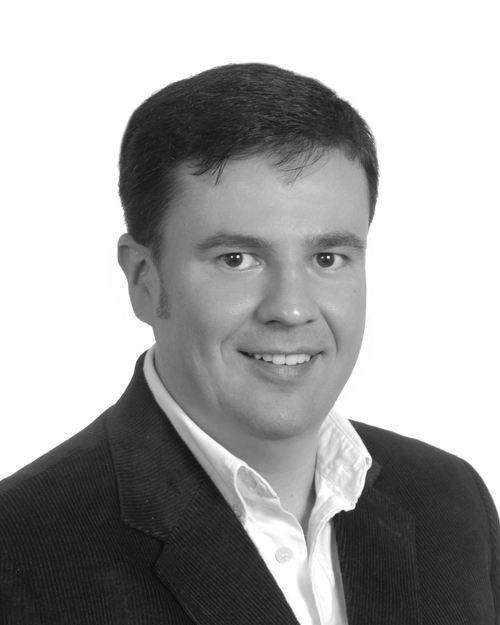}}]{Miguel Ortega-Vazquez}
	(S'97-M'06-SM'14) is a Senior Technical Leader at the Electric Power Research Institute in Palo Alto, CA.  His current research interests include power system operation, power system security, power system economics, integration of renewable energy sources, and the smart grid.  Before joining EPRI, he was Assistant Professor at the University of Washington in Seattle, WA; and before that, Assistant Professor at the Chalmers University of Technology, Sweden.  He holds a Ph.D. from The University of Manchester; M.Sc. from the Universidad Aut\'{o}noma de Nuevo Le\'{o}n, Mexico; and an Electric Engineering degree from the Instituto Tecnol\'{o}gico de Morelia, Mexico.
	\end{IEEEbiography}

\end{document}